\documentclass[12pt]{article}%
\usepackage{amsfonts}
\usepackage{sw20bams}
\usepackage{amsmath}
\usepackage{amssymb}
\usepackage{graphicx}%
\providecommand{\U}[1]{\protect\rule{.1in}{.1in}}
\begin{document}

\title{D/M/1 Queue: Policies and Control}
\author{Steven Finch}
\date{October 27, 2022}
\maketitle

\begin{abstract}
Equilibrium G/M/1-FIFO waiting times are exponentially distributed, as first
proved by Smith (1953). \ For other client-sorting policies, such generality
is not feasible. Assume that interarrival times are constant. \ Symbolics for
the D/M/1-LIFO\ density are completely known; numerics for D/M/1-SIRO arise
via an unpublished recursion due to Burke (1967). \ Consider a weighted sum of
two costs, one from keeping clients waiting for treatment and the other from
having the server idle. \ With this in mind, what is the optimal interarrival
time and how does this depend on the choice of policy?

\end{abstract}

\footnotetext{Copyright \copyright \ 2022 by Steven R. Finch. All rights
reserved.}In an M/G/1 system, interarrival times are exponentially distributed
with mean $1/\lambda$ and service lengths are arbitrary with mean $1/\mu$.
\ In a G/M/1 system, it is service lengths that are exponentially distributed
with mean $1/\mu$, while interarrival times are arbitrary with mean
$1/\lambda$. \ All such intervals are taken to be independent; client waiting
space is unlimited. \ The two systems display a kind of subtle symmetry, never
overt. \ Let $\rho=\lambda/\mu$ denote the traffic intensity. \ Our focus in
this paper is D/M/1, for which interarrival times are constant $a=1/\lambda$.

We discussed first-in-first-out M/D/1 in \cite{Fi1-que3}; both
last-in-first-out M/D/1 and serve-in-random-order M/D/1 were covered in
\cite{Fi2-que3}. \ The two policies FIFO and LIFO\ are also known as FCFS
(first-come-first-serve) and LCFS (last-come-first-serve). \ The policy
SIRO\ is also variously known as ROS (random-order-of-service) and RSS (random-selection-for-service).

Given D/M/1, the functional equation%
\[
\zeta(s)=\exp\left[  -a\left(  \mu+s-\mu\,\zeta(s)\right)  \right]
\]
is important. \ This is a consequence of $\exp(-a\,s)$ being the Laplace
transform of $\delta(x-a)$, where $\delta(x)$ is the Dirac delta. \ An
analogous equation holds for non-constant interarrival times, depending on the
probability density, but this would be more complicated. \ We have%
\[
\zeta(s)=-\frac{1}{a\,\mu}\,\omega\left[  -a\,\mu\exp\left(  -a\,\mu
-a\,s\right)  \right]  =-\rho\,\omega\left[  -\frac{1}{\rho}\exp\left(
-\frac{1}{\rho}-a\,s\right)  \right]
\]
where $\omega(s)$ is the principal branch of the Lambert omega:
\[%
\begin{array}
[c]{ccccc}%
\omega(s)e^{\omega(s)}=s, &  & -1\leq\omega(x)\in\mathbb{R}\text{ \ }%
\forall\text{ }x\geq-1/e, &  & \exists\text{ branch cut for }x<-1/e.
\end{array}
\]
As an example of aforementioned symmetry, the expression for $\zeta(s)$ is
identical to that for $\Theta(s)$ corresponding to M/D/1 in \cite{Fi2-que3}
except $\rho$ is everywhere replaced by $1/\rho$. \ A\ special value%
\[
\zeta_{0}=\zeta(0)=-\frac{1}{a\,\mu}\,\omega\left(  -a\,\mu\,e^{-a\,\mu
}\right)
\]
will appear throughout. \ For instance, if $\lambda=2$ and $\mu=3$, then
$\zeta_{0}=0.41718835...$. \ Upon differentiation, we have%
\[%
\begin{array}
[c]{ccc}%
\omega^{\prime}(s)=\dfrac{\omega(s)}{s\left(  1+\omega(s)\right)  }, &  &
\zeta^{\prime}(s)=-a\,\omega^{\prime}\left[  -a\,\mu\exp\left(  -a\,\mu
-a\,s\right)  \right]  \exp\left(  -a\,\mu-a\,s\right)  .
\end{array}
\]
Another value%
\begin{align*}
\zeta_{0}^{\prime}  &  =\zeta^{\prime}(0)=-a\,\frac{\omega\left(
-a\,\mu\,e^{-a\,\mu}\right)  }{\left(  -a\,\mu\,e^{-a\,\mu}\right)  \left[
1+\omega\left(  -a\,\mu\,e^{-a\,\mu}\right)  \right]  }\,e^{-a\,\mu}\\
&  =-a\,\frac{-a\,\mu\,\zeta_{0}}{\left(  -a\,\mu\right)  \left(
1-a\,\mu\,\zeta_{0}\right)  }=-\frac{a\,\zeta_{0}}{1-a\,\mu\,\zeta_{0}}%
\end{align*}
is also needed; if $\lambda=2$ and $\mu=3$, then $\zeta_{0}^{\prime
}=-0.55741433...$.

\section{FIFO}

Let $W_{\text{que}}$ denote the waiting time in the queue (prior to service).
\ Under equilibrium (steady-state) conditions, the probability density
function $f(x)$ of $W_{\text{que}}$ has Laplace transform \cite{Smi-que3,
Pol-que3, Ken-que3, Wis-que3}
\[
F(s)=\lim_{\varepsilon\rightarrow0^{+}}%
{\displaystyle\int\limits_{-\varepsilon}^{\infty}}
\exp(-s\,x)f(x)dx=1-\zeta_{0}+\zeta_{0}\,\frac{\mu(1-\zeta_{0})}{s+\mu
(1-\zeta_{0})}=F_{\text{alt}}(s)+1-\zeta_{0}%
\]
and initial value \cite{Sch-que3}%
\[
f(0^{+})=\lim_{s\rightarrow1\cdot\infty}s\,F_{\text{alt}}(s)=\mu\,\zeta
_{0}\,(1-\zeta_{0}).
\]
Consequently%
\[
f(x)=(1-\zeta_{0})\delta(x)+\mu\,\zeta_{0}\,(1-\zeta_{0})\exp\left(
-\mu(1-\zeta_{0})x\right)  .
\]
In fact, exponentiality holds more generally for non-constant interarrival
times, proved by Smith \cite{Smi-que3}. Moments are%
\[%
\begin{array}
[c]{ccc}%
\text{mean}=-F^{\prime}(0)=\dfrac{\zeta_{0}}{\mu(1-\zeta_{0})}, &  &
\text{variance}=F^{\prime\prime}(0)-F^{\prime}(0)^{2}=\zeta_{0}\,\dfrac
{2-\zeta_{0}}{\mu^{2}(1-\zeta_{0})^{2}}%
\end{array}
\]
giving $0.23860673...$ and $0.21600433...$ respectively when $\{\lambda
,\mu\}=\{2,3\}$. \ If sampling is restricted only to $W_{\text{que}}>0$, then
\cite{Tak-que3}%
\[%
\begin{array}
[c]{ccc}%
\text{mean}_{>0}=\dfrac{1}{\mu(1-\zeta_{0})}, &  & \text{variance}_{>0}%
=\dfrac{1}{\mu^{2}(1-\zeta_{0})^{2}}%
\end{array}
\]
giving $0.57194007...$ and $0.32711544...$ respectively.

Let $L_{\text{sys}}$ denote the number of patients in the system (both queue
and service). \ Under equilibrium, with $\{\lambda,\mu\}=\{2,3\}$, we have
\cite{Bha-que3}%
\[%
\begin{array}
[c]{ccc}%
\tilde{f}(\ell)=\mathbb{P}\left\{  L_{\text{sys}}=\ell\right\}  =(1-\zeta
_{0})\zeta_{0}^{\ell}, &  & \ell=0,1,2,3,\ldots;
\end{array}
\]%
\[%
\begin{array}
[c]{ccc}%
\text{mean}=\dfrac{\zeta_{0}}{1-\zeta_{0}}, &  & \text{variance}=\dfrac
{\zeta_{0}}{\left(  1-\zeta_{0}\right)  ^{2}}%
\end{array}
\]
giving $0.71582021...$ and $1.22821879...$ respectively. \ Geometricity holds
more generally for non-constant interarrival times. \ It is remarkable that
classical distributions occur within G/M/1 universally but not within even
M/D/1 specifically.

\section{LIFO}

The probability density function $f(x)$ of $W_{\text{que}}$ has Laplace
transform
\[
F(s)=\lim_{\varepsilon\rightarrow0^{+}}%
{\displaystyle\int\limits_{-\varepsilon}^{\infty}}
\exp(-s\,x)f(x)dx=1-\zeta_{0}+\zeta_{0}\frac{\mu-\mu\,\zeta(s)}{s+\mu
-\mu\,\zeta(s)}=F_{\text{alt}}(s)+1-\zeta_{0}%
\]
and the inverse Laplace transform of $\zeta(s)$ is%
\[
\theta(x)=%
{\displaystyle\sum\limits_{k=1}^{\infty}}
e^{-k/\rho}\dfrac{(k/\rho)^{k-1}}{k!}\delta(x-a\,k).
\]
With regard to symmetry and $F(s)$, we see that $\{\zeta_{0},\mu,\zeta\}$ play
the roles of $\{\rho,\lambda,\Theta\}$ in \cite{Fi2-que3}, but an extra factor
$\zeta_{0}$ is also present, i.e., the correspondence is not perfect. \ From%
\[
(1-\zeta_{0})s+\mu\,\left[  1-\zeta(s)\right]  (1-\zeta_{0}+\zeta
_{0})=s\,F(s)+\mu\,F(s)\left[  1-\zeta(s)\right]
\]
we have%
\[
(1-\zeta_{0})s+\mu\left(  1-F(s)\right)  \left[  1-\zeta(s)\right]  =s\,F(s),
\]
i.e.,%
\[
F(s)=1-\zeta_{0}+\mu\left(  1-F(s)\right)  \left[  \frac{1}{s}-\frac{\zeta
(s)}{s}\right]
\]
hence%
\begin{align*}
f(x)  &  =(1-\zeta_{0})\delta(x)+\kappa+\mu%
{\displaystyle\int\limits_{0}^{x}}
\left(  \delta(t)-f(t)\right)  \left[  1-%
{\displaystyle\int\limits_{0}^{x-t}}
\theta(u)du\right]  dt\\
&  =(1-\zeta_{0})\delta(x)+\mu\,\zeta_{0}+\mu\left[  1-%
{\displaystyle\int\limits_{0}^{x}}
\theta(u)du\right]  -\mu%
{\displaystyle\int\limits_{0}^{x}}
f(t)\left[  1-%
{\displaystyle\int\limits_{0}^{x-t}}
\theta(u)du\right]  dt.
\end{align*}
The indicated condition $\kappa=\mu\,\zeta_{0}$ is true by the initial value
theorem \cite{Sch-que3}:%
\[
\lim_{\varepsilon\rightarrow0^{+}}f(\varepsilon)=\lim\limits_{s\rightarrow
1\cdot\infty}s\,F_{\text{alt}}(s).
\]
Differentiating, we obtain
\begin{align*}
f^{\prime}(x)  &  =\mu\left[  0-\theta(x)\right]  -\mu f(x)\left[  1-0\right]
-\mu%
{\displaystyle\int\limits_{0}^{x}}
f(t)\left[  0-\theta(x-t)\right]  dt\\
&  =-\mu\theta(x)-\mu f(x)+\mu%
{\displaystyle\int\limits_{0}^{x}}
f(t)\theta(x-t)dt\\
&  =-\mu\theta(x)-\mu f(x)+\mu%
{\displaystyle\int\limits_{0}^{x}}
f(t)%
{\displaystyle\sum\limits_{k=1}^{\infty}}
e^{-k/\rho}\dfrac{(k/\rho)^{k-1}}{k!}\delta(x-t-a\,k)dt\\
&  =-\mu\theta(x)-\mu f(x)+\mu%
{\displaystyle\sum\limits_{k=1}^{\infty}}
e^{-k/\rho}\dfrac{(k/\rho)^{k-1}}{k!}%
{\displaystyle\int\limits_{0}^{x}}
f(t)\delta(x-t-a\,k)dt\\
&  =-\mu\theta(x)-\mu f(x)+\mu%
{\displaystyle\sum\limits_{k=1}^{\infty}}
e^{-k/\rho}\dfrac{(k/\rho)^{k-1}}{k!}f(x-a\,k).
\end{align*}
For $0<x<a$,%
\[%
\begin{array}
[c]{ccc}%
f^{\prime}(x)=-\mu\,f(x), &  & f(0^{+})=\mu\,\zeta_{0}%
\end{array}
\]
implies%
\[
f(x)=\mu\,\zeta_{0}\,e^{-\mu\,x}.
\]
Note that $\lim_{\varepsilon\rightarrow0^{+}}f(a\,k+\varepsilon)=0$ for each
$k\geq1$ because, if a client arrives at the same moment the server becomes
available, the client is taken immediately (by LIFO) and there is no waiting.
Note also $1/\rho=\mu/\lambda=\mu\,a$. \ For $a<x<2a$,%
\begin{align*}
f^{\prime}(x)  &  =-\mu\,f(x)+\mu\,e^{-1/\rho}\cdot\mu\,\zeta_{0}%
\,e^{-\mu(x-a)}\\
&  =-\mu\,f(x)+\mu^{2}\zeta_{0}\,e^{-\mu\,x}%
\end{align*}
coupled with $f(a^{+})=0$ implies%
\[
f(x)=\mu^{2}\zeta_{0}(x-a)e^{-\mu\,x}.
\]
For $2a<x<3a$,%
\begin{align*}
f^{\prime}(x)  &  =-\mu\,f(x)+\mu\,e^{-1/\rho}\cdot\mu^{2}\zeta_{0}%
(x-2a)e^{-\mu(x-a)}+\mu\,e^{-2/\rho}\frac{2/\rho}{2!}\cdot\mu\,\zeta
_{0}e^{-\mu(x-2a)}\\
&  =-\mu\,f(x)+\mu^{3}\zeta_{0}(x-2a)e^{-\mu\,x}+\left(  \mu^{2}\,\zeta
_{0}\right)  (\mu\,a)e^{-\mu\,x}\\
&  =-\mu\,f(x)+\mu^{3}\zeta_{0}(x-a)e^{-\mu\,x}%
\end{align*}
coupled with $f(2a^{+})=0$ implies%
\[
f(x)=\frac{1}{2}\mu^{3}\zeta_{0}\,x(x-2a)e^{-\mu\,x}.
\]
More generally, for $k\,a<x<(k+1)a$, we obtain%
\[
f(x)=\frac{1}{k!}\mu^{k+1}\zeta_{0}\,x^{k-1}(x-k\,a)e^{-\mu\,x}%
\]
and thus the waiting time density for LIFO\ is completely understood.
\ Wishart\ \cite{Wis-que3} evidently holds priority in discovering this
formula, building upon work by Conolly \cite{Con1-que3}. Stitching the
fragments together gives the LIFO density function pictured in Figure 1, for
parameter values $\lambda=2$ and $\mu=3$; hence $\rho=2/3$ and $a=1/2$. \ 

Moments of $W_{\text{que}}$ for LIFO\ are \cite{Wis-que3}\
\[%
\begin{array}
[c]{ccc}%
\text{mean}=-F^{\prime}(0)=\dfrac{\zeta_{0}}{\mu(1-\zeta_{0})}, &  &
\text{variance}=F^{\prime\prime}(0)-F^{\prime}(0)^{2}=\zeta_{0}\,\dfrac
{2-\zeta_{0}-2\mu\,\zeta_{0}^{\prime}}{\mu^{2}(1-\zeta_{0})^{2}}%
\end{array}
\]
giving $0.23860673...$ and $0.67242217...$ respectively. \ The mean of
$W_{\text{que}}$ for FIFO is the same as that for LIFO; the variance for FIFO
is smaller. \ If sampling is restricted only to $W_{\text{que}}>0$, then
\cite{Tak-que3, Pra-que3, HC-que3, CH-que3}%
\[%
\begin{array}
[c]{ccc}%
\text{mean}_{>0}=\dfrac{1}{\mu(1-\zeta_{0})}, &  & \text{variance}_{>0}%
=\dfrac{1-2\mu\,\zeta_{0}^{\prime}}{\mu^{2}(1-\zeta_{0})^{2}}%
\end{array}
\]
giving $0.57194007...$ and $1.42114846...$ respectively. \ The variance
expression reported in \cite{Con2-que3} contains an apparent error.

\section{SIRO}

The probability density function $f(x)$ of $W_{\text{que}}$ has Laplace
transform \cite{Tak-que3}%
\[
F(s)=1-\zeta_{0}+\zeta_{0}\,\Phi(s)=F_{\text{alt}}(s)+1-\zeta_{0}%
\]
where%
\[
\Phi(s)=B(s,\zeta_{0})-%
{\displaystyle\int\limits_{\zeta(s)}^{\zeta_{0}}}
\exp\left(  -%
{\displaystyle\int\limits_{u}^{\zeta_{0}}}
\frac{dv}{v-e^{-a\left(  \mu+s-\mu\,v\right)  }}\right)  \frac{\partial
B}{\partial u}(s,u)du,
\]%
\[
B(s,z)=\frac{\mu(1-\zeta_{0})}{1-z}\,\frac{1-\exp\left[  -a\left(  s+\mu
-\mu\,z\right)  \right]  }{s+\mu-\mu\,z}.
\]
The integral underlying $\Phi(s)$ is intractable; our symbolic approach for
FIFO\ \&\ LIFO seems inapplicable for SIRO.

We therefore turn to a numeric approach. \ An unpublished memorandum written
in 1967 by Burke (the same author as of \cite{Bur-que3}) has regrettably been
lost, although summaries are found in \cite{CC-que3, Coo-que3}. \ Rosenlund
\cite{Ros-que3} provided an especially clear algorithm for D/M/1 to follow.
\ Since our interest is in densities, we differentiate his initial expression
with respect to $x$, i.e.,%
\[
\frac{d}{dx}\left(  x^{j+1-r}e^{-x}\right)  =(j+1-r-x)x^{j-r}e^{-x}.
\]
Define recursively%
\[%
\begin{array}
[c]{ccc}%
h_{j,0}(x)=\,%
{\displaystyle\sum\limits_{r=1}^{j+1}}
\dfrac{r}{j+1}\,\dfrac{(j+1-r-x)x^{j-r}}{(j+1-r)!}\,e^{-x}, &  &
j=0,1,2,\ldots;
\end{array}
\]%
\[%
\begin{array}
[c]{ccc}%
h_{j,k}(x)=\,%
{\displaystyle\sum\limits_{r=1}^{j+1}}
\dfrac{r}{j+1}\,\dfrac{(1/\rho)^{j+1-r}}{(j+1-r)!}\,e^{-1/\rho}\,h_{r,k-1}%
(x), &  & j=0,1,2,\ldots\text{ and }k=1,2,3,\ldots.
\end{array}
\medskip
\]
We consequently have%
\[
f(x)=(1-\zeta_{0})\delta(x)+\zeta_{0}\,g(x)
\]
where%
\[%
\begin{array}
[c]{ccc}%
g(x)=-\mu(1-\zeta_{0})%
{\displaystyle\sum\limits_{j=0}^{\infty}}
\,\zeta_{0}^{j}\,h_{j,\left\lfloor \lambda\,x\right\rfloor }\left(
\mu\,x-\dfrac{\left\lfloor \lambda\,x\right\rfloor }{\rho}\right)  , &  &
x\geq0.
\end{array}
\]
For example, if $0<x<a$, then
\[
f(x)=\mu(1-\zeta_{0})%
{\displaystyle\int\limits_{1-\zeta_{0}}^{1}}
\frac{e^{-\mu\,x\,t}}{t}dt=\mu(1-\zeta_{0})\left[  E\left(  \mu(1-\zeta
_{0})x\right)  -E\left(  \mu\,x\right)  \right]
\]
where $E(x)=-\operatorname{Ei}(-x)$ is the exponential integral.\ \ This
corresponds to the leftmost curvilinear arc in Figure 2, surmounting the
interval $[0,\frac{1}{2}]$. \ Verification that the Laplace transform of
$f(x)$ is equal to $F_{\text{alt}}(s)$ remains open.

It is known (by other techniques) that the mean of $W_{\text{que}}$ for SIRO
is the same as that for FIFO and LIFO; the corresponding variance is between
the two extremes \cite{Tak-que3, Ros-que3}:
\[
\zeta_{0}\,\dfrac{4-2\zeta_{0}-4\mu\,\zeta_{0}^{\prime}+\mu\,\zeta_{0}%
\,\zeta_{0}^{\prime}}{\mu^{2}(1-\zeta_{0})^{2}\left(  2-\mu\,\zeta_{0}%
^{\prime}\right)  }%
\]
giving $0.34029290...$. \ If sampling is restricted only to $W_{\text{que}}%
>0$, then the variance is
\[
\dfrac{2-3\mu\,\zeta_{0}^{\prime}}{\mu^{2}(1-\zeta_{0})^{2}\left(
2-\mu\,\zeta_{0}^{\prime}\right)  }%
\]
giving $0.62503500...$.

\section{Idle Period}

We are concerned here with successive periods of server activity and
inactivity. \ The left-hand subinterval of $k\,a\leq x<(k+1)a$ is busy (since
a new client has just arrived) and its right-hand complement is idle. \ It is
possible that the idle period is empty. \ Jansson \cite{Jan-que3} proved that,
under FIFO and equilibrium, the idle period length has probability density
function%
\[%
\begin{array}
[c]{ccc}%
\zeta_{0}\,\delta(x)+\mu\,\zeta_{0}\,(1-\zeta_{0})\exp\left(  \mu(1-\zeta
_{0})x\right)  , &  & 0\leq x<a.
\end{array}
\]
Moments are%
\[%
\begin{array}
[c]{ccc}%
\text{mean}=\dfrac{1}{\lambda}-\dfrac{1}{\mu}=(1-\rho)a, &  & \text{variance}%
=\dfrac{1+\zeta_{0}-2a\,\mu\,\zeta_{0}}{\mu^{2}(1-\zeta_{0})}%
\end{array}
\]
giving $1/6$ and $0.03157553...$ respectively when $\{\lambda,\mu\}=\{2,3\}$.
\ The analysis of a busy period is more complicated, in part because it may
span multiple adjacent intervals $[k\,a,\,(k+1)a)$, but this issue is not
pertinent for our study here.

Each client is associated with both a waiting time $\in\lbrack0,\infty)$ and
an idle period length $\in\lbrack0,a)$. \ An expression for the\ bivariate
density is available \cite{Jan-que3}. \ We report merely the cross-covariance%

\[
\dfrac{(1-\zeta_{0})e^{-a\,\mu}}{\mu^{2}\zeta_{0}}+\dfrac{a\,\zeta_{0}}%
{\mu(1-\zeta_{0})}-\dfrac{1}{\mu^{2}(1-\zeta_{0})}%
\]
and cross-correlation $-0.44448913...$ when $\{\lambda,\mu\}=\{2,3\}$.
\ Again, the proof is valid under FIFO and equilibrium. \ What is remarkable
is that these results (marginal density and joint moments) appear via
simulation to be the same under LIFO\ and\ SIRO as well. \ Likewise, the
distribution of $L_{\text{sys}}$ (what we called $\tilde{f}$ in Section 1)
seems to be invariant upon change in policy. \ Justification would be good to
see someday.

\section{Minimal Cost}

The expression \textquotedblleft queue control\textquotedblright\ may seem
redundant because queues are themselves a method of control \cite{Dai-que3}.
\ They exist to accommodate client demands on a service provider.
\ A\ control, however, exists to ensure that costs remain sustainable. \ We
wish to minimize cost as a function of $a=1/\lambda$, for fixed $\mu$, where
cost is a $c$-weighted sum of the mean idle period and the mean waiting time
\cite{Jan-que3}:
\[
C=(1-c)\left(  a-\frac{1}{\mu}\right)  +c\,\dfrac{\zeta_{0}}{\mu(1-\zeta_{0}%
)}.
\]
The derivative of $\zeta_{0}$ with respect to $a$ will be written as
$\zeta_{0}^{\prime}$, which should not be confused with our earlier usage of
the same symbol (the derivative of $\zeta$ with respect to $s$, evaluated at
$0$). \ From%
\[
\zeta_{0}=-\frac{1}{a\,\mu}\,\omega\left(  -a\,\mu\,e^{-a\,\mu}\right)
\]
we deduce%
\begin{align*}
\left(  \mu\,\zeta_{0}\right)  ^{\prime}  &  =\frac{1}{a^{2}}\,\omega\left(
-a\,\mu\,e^{-a\,\mu}\right)  -\frac{1}{a}\,\omega^{\prime}\left(
-a\,\mu\,e^{-a\,\mu}\right)  \left(  -\mu\,e^{-a\,\mu}+a\,\mu^{2}e^{-a\,\mu
}\right) \\
&  =\frac{1}{a^{2}}\,\omega\left(  -a\,\mu\,e^{-a\,\mu}\right)  -\frac{1}%
{a}\,\frac{\omega\left(  -a\,\mu\,e^{-a\,\mu}\right)  }{\left(  -a\,\mu
\,e^{-a\,\mu}\right)  \left[  1+\omega\left(  -a\,\mu\,e^{-a\,\mu}\right)
\right]  }\left(  -\mu\right)  (1-a\,\mu)e^{-a\,\mu}\\
&  =\frac{\omega\left(  -a\,\mu\,e^{-a\,\mu}\right)  }{a^{2}}\,\left[
1-\frac{1-a\,\mu}{1+\omega\left(  -a\,\mu\,e^{-a\,\mu}\right)  }\right]
=\frac{\omega\left(  -a\,\mu\,e^{-a\,\mu}\right)  }{a^{2}}\,\frac
{a\,\mu+\omega\left(  -a\,\mu\,e^{-a\,\mu}\right)  }{1+\omega\left(
-a\,\mu\,e^{-a\,\mu}\right)  }\\
&  =-\frac{\mu\,\zeta_{0}}{a}\,\frac{a\,\mu-a\,\mu\,\zeta_{0}}{1-a\,\mu
\,\zeta_{0}}=-\frac{\mu^{2}\zeta_{0}\left(  1-\zeta_{0}\right)  }%
{1-a\,\mu\,\zeta_{0}}%
\end{align*}
thus%
\begin{align*}
\left(  \dfrac{\zeta_{0}}{1-\zeta_{0}}\right)  ^{\prime}  &  =\dfrac{\zeta
_{0}^{\prime}}{1-\zeta_{0}}-\dfrac{\zeta_{0}}{(1-\zeta_{0})^{2}}\left(
-\zeta_{0}^{\prime}\right) \\
&  =\left[  \dfrac{1}{1-\zeta_{0}}+\dfrac{\zeta_{0}}{(1-\zeta_{0})^{2}%
}\right]  \,\left(  -\frac{\mu\,\zeta_{0}\left(  1-\zeta_{0}\right)
}{1-a\,\mu\,\zeta_{0}}\right) \\
&  =-\dfrac{1-\zeta_{0}+\zeta_{0}}{(1-\zeta_{0})^{2}}\,\frac{\mu\,\zeta
_{0}\left(  1-\zeta_{0}\right)  }{1-a\,\mu\,\zeta_{0}}=\frac{-\mu\,\zeta_{0}%
}{\left(  1-\zeta_{0}\right)  \left(  1-a\,\mu\,\zeta_{0}\right)  }%
\end{align*}
thus%
\[
C^{\prime}=(1-c)-c\,\frac{\zeta_{0}}{\left(  1-\zeta_{0}\right)  \left(
1-a\,\mu\,\zeta_{0}\right)  }=0
\]
when%
\[
\frac{c}{1-c}=\frac{\left(  1-\zeta_{0}\right)  \left(  1-a\,\mu\,\zeta
_{0}\right)  }{\zeta_{0}}=\frac{1-\zeta_{0}-\zeta_{0}\cdot a\,\mu(1-\zeta
_{0})}{\zeta_{0}}.
\]
It is additionally required \cite{Jan-que3} that $0<\zeta_{0}<1$. \ From%
\[%
\begin{array}
[c]{ccccc}%
\zeta_{0}=\exp\left(  -a\,\mu(1-\zeta_{0})\right)  , &  & \text{i.e.,} &  &
a\,\mu(1-\zeta_{0})=-\ln(\zeta_{0})
\end{array}
\]
we obtain%
\[
\frac{c}{1-c}=\frac{1-\zeta_{0}+\zeta_{0}\cdot\ln(\zeta_{0})}{\zeta_{0}}%
\]
hence%
\[
\zeta_{0}=-\frac{1}{\bar{\omega}\left[  -\exp\left(  -\dfrac{1}{1-c}\right)
\right]  }%
\]
where $\bar{\omega}(s)$ is \textquotedblleft the\textquotedblright\ secondary
branch of the Lambert omega:
\[%
\begin{array}
[c]{ccccc}%
\bar{\omega}(s)e^{\bar{\omega}(s)}=s, &  & -1\geq\bar{\omega}(x)\in
\mathbb{R}\text{ \ }\forall\text{ }x\in\lbrack-1/e,0), &  & \exists\text{
branch cut for }x\leq0.
\end{array}
\]
For example, if $\mu=3$ and $c=1/2$, then $\zeta_{0}=0.31784443...$ and
$a=0.56008398...$. \ In words, if mean client waiting times are weighted the
same as mean server idle periods, i.e., $c=1/2$, then in terms of cost, the
interarrival time $a=0.50$ is far from optimal, but $a=0.56$ is close.

If server idle periods are weighted more heavily than client waiting times,
e.g., $c=1/5$, then $a=0.44983251...$. \ If instead client waiting times are
weighted more heavily than server idle periods, e.g., $c=4/5$, then
$a=0.75436304...$. \ This is consistent with intuition. \ Compressed
interarrival times lead to less idleness but longer waits; expansive
interarrival times lead to shorter waits but more idleness. \ Balancing these
conflicting priorities makes life interesting.

To clarify: there exist countably infinite branches of the Lambert omega, but
only two ($\omega$ and $\bar{\omega}$) that assume real values on $[-1/e,0)$,
one increasing and the other decreasing. \ All other branches are
complex-valued with nonzero imaginary parts. \ Our notation $\bar{\omega}$ is
unorthodox, as is referring to $\bar{\omega}$ as \textquotedblleft
the\textquotedblright\ secondary branch. \ In Mathematica, the function
\texttt{ProductLog[k,x]} gives $\omega(x)$ \&\ $\bar{\omega}(x)$ for $k=0$ \&
$k=-1$, respectively. \ Alternative notation $\omega_{+}$ and $\omega_{-}$,
proposed somewhat by \cite{NIST-que3}, is intended to suggest
\textquotedblleft upper branch\textquotedblright\ and \textquotedblleft lower
branch\textquotedblright.

We have omitted discussion of the variance of $C$. \ From the aforementioned
joint distribution of idle period and waiting time \cite{Jan-que3}, it would
be possible to minimize cost as a function of $a$, for fixed $\mu$, where cost
is the median of a $c$-weighted sum of idle period and waiting time. \ Solving
this revised optimization problem could be advantageous because the median is
more a robust estimator of centrality than the mean. \ We wonder too about the
proper choice of $0<c<1$ and whether a sum (rather than a product, say) is
necessarily best. \ More recent work appears in \cite{X1-que3, X2-que3,
X3-que3, X4-que3, X5-que3}. \ Processes with constant interarrival times and
exponential server queues are fundamental, as proved in \cite{Haj-que3}.%

\begin{figure}
[ptb]
\begin{center}
\includegraphics[
height=4.0568in,
width=6.2379in
]%
{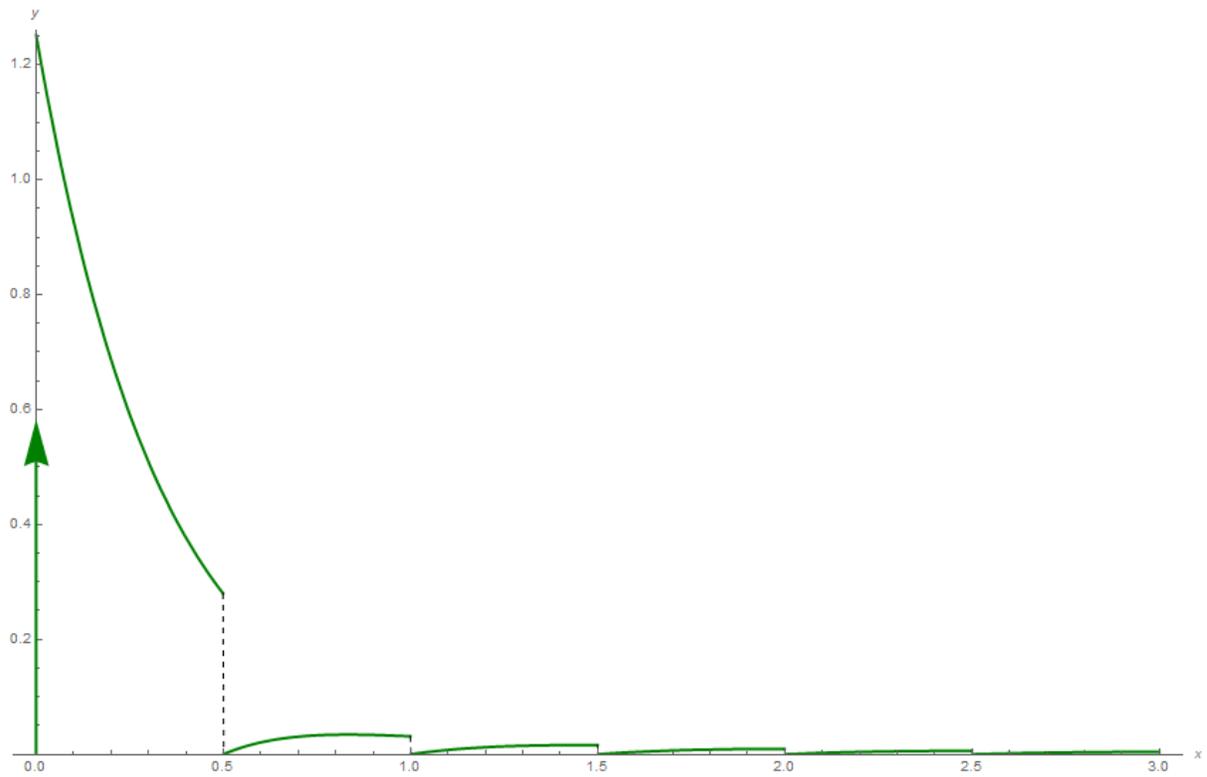}%
\caption{Waiting time density plot for Deterministic$[\frac{1}{2}]$ arrivals,
last-in-first-out exponential service.}%
\end{center}
\end{figure}
\begin{figure}
[ptb]
\begin{center}
\includegraphics[
height=4.2142in,
width=6.5008in
]%
{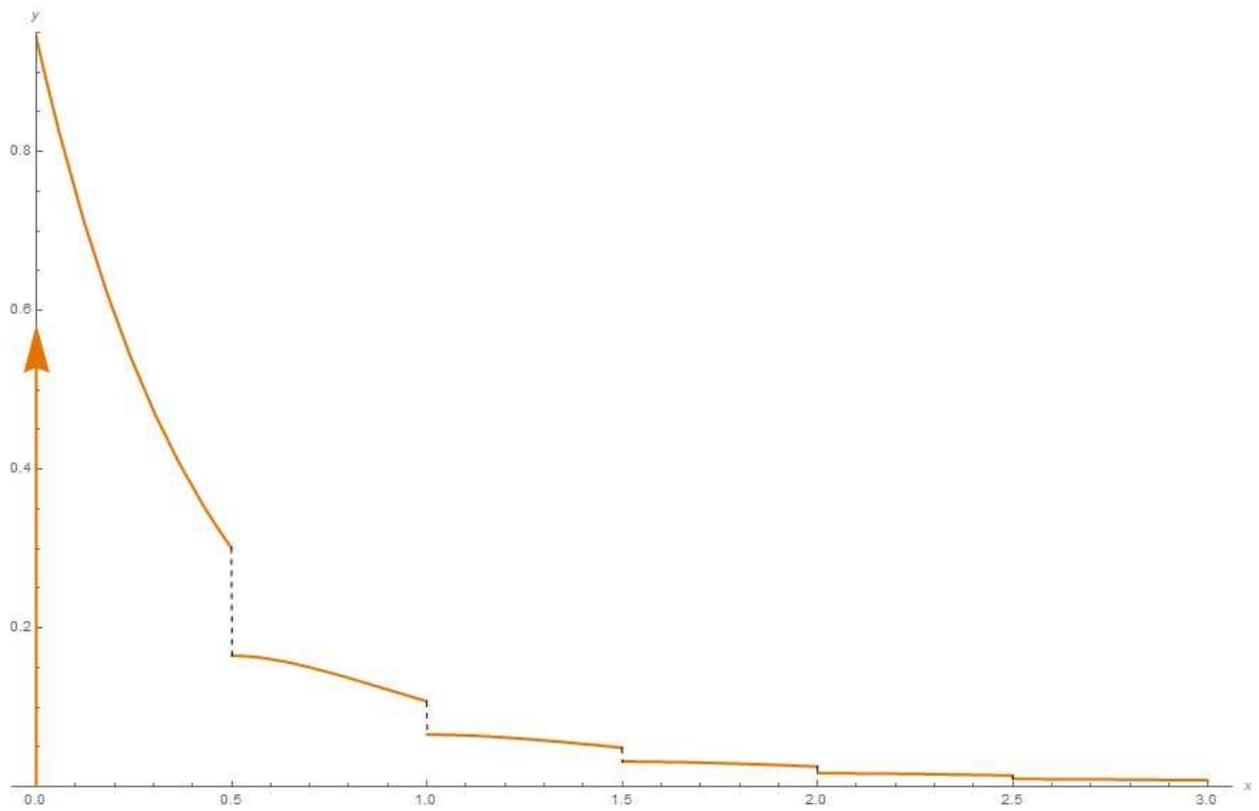}%
\caption{Waiting time density plot for Deterministic$[\frac{1}{2}]$ arrivals,
exponential serve-in-random order policy.}%
\end{center}
\end{figure}

\section{Addendum}

With $\zeta_{0}$ as before (constant $\approx0.417$ for $\lambda=2=1/a$,
$\mu=3$), define \cite{Ken-que3, Tak-que3, Tkc-que3, All-que3}%
\[
\Delta_{n}(x)=\left\{  \frac{1}{1-\zeta_{0}}+%
{\displaystyle\sum\limits_{j=1}^{n}}
\dbinom{n}{j}\frac{Q_{j}(x)}{1-e^{-a\,j\,x}}\,\frac{n\left(  1-e^{-a\,j\,x}%
\right)  -j}{n\left(  1-\zeta_{0}\right)  -j}\right\}  ^{-1},
\]%
\[%
\begin{array}
[c]{ccc}%
m_{n}=\dfrac{\Delta_{n}(\mu/n)}{\mu\left(  1-\zeta_{0}\right)  ^{2}}, &  &
p_{n}=1-\dfrac{\Delta_{n}(\mu/n)}{1-\zeta_{0}}%
\end{array}
\]
where $n$ is a positive integer and%
\[
Q_{j}(x)=%
{\displaystyle\prod\limits_{i=1}^{j}}
\,\frac{1-e^{-a\,i\,x}}{e^{-a\,i\,x}}.
\]
For example, $\Delta_{1}=\left(  1-\zeta_{0}\right)  \zeta_{0}$, $m_{1}%
=\zeta_{0}/\left(  \mu\left(  1-\zeta_{0}\right)  \right)  \approx0.239$ and
$p_{1}=1-\zeta_{0}\approx0.583$. \ More generally, $m_{n}$ is the expected
waiting time in a D/M/$n$ queue with $n$ slow servers (more precisely, each
server working with rate only $\mu/n$ when busy) and $p_{n}$ is the
probability of zero wait. \ With $\{\lambda,\mu\}=\{2,3\}$, we have%
\[%
\begin{array}
[c]{ccc}%
m_{2}=0.16901950..., &  & p_{2}=0.70448039...;
\end{array}
\]%
\[%
\begin{array}
[c]{ccc}%
m_{3}=0.12647170..., &  & p_{3}=0.77887245...;
\end{array}
\]%
\[%
\begin{array}
[c]{ccc}%
m_{4}=0.09744181..., &  & p_{4}=0.82962932...;
\end{array}
\]%
\[%
\begin{array}
[c]{ccc}%
m_{5}=0.07648770..., &  & p_{5}=0.89364299...;
\end{array}
\]
i.e., $n$ slow servers outperform one fast server, relative to average waiting
time. \ The sum $S$ of idle periods over all servers would however be
potentially significant; the mean of $S$ would be crucial in minimizing total
cost as a function of $a$, for fixed $\mu$.

\section{Acknowledgements}

Stig Rosenlund and Robert Cooper were so kind in answering several of my
questions. I am grateful to innumerable software developers. \ Mathematica
routines NDSolve for delay-differential equations and InverseLaplaceTransform
(for Mma version $\geq12.2$) assisted in numerically confirming many results.
\ R steadfastly remains my favorite statistical programming language.

\end{document}